\def\thetitle{Completeness of the ring of polynomials}
\def\thedate{13 December 2013}

\input amstex
\documentstyle{amsppt}
\magnification=\magstep1

%
%

\def\cref#1{#1}
\def\pref#1{(#1)}

\def\sectno#1{\relax}

\chardef\atcode=\catcode`\@
\catcode`\@=11
\def\varprojlim{\mathop{\vtop{\ialign{##\crcr
 \hfil\rm lim\hfil\crcr
 \noalign{\nointerlineskip\vskip1pt }\leftarrowfill\crcr
 \noalign{\nointerlineskip\kern-\ex@}\crcr}}}}
\catcode`\@=\atcode

\def\Hom{\operatorname{Hom}}
\def\Ext{\operatorname{Ext}}

\document
\topmatter
   \date \thedate \enddate
 \title  \thetitle  \endtitle
 \author Anders Thorup
 \endauthor
\affil 
Department of Mathematics, University of Copenhagen
 \endaffil 
\address  Universitetsparken 5, DK--2100 Copenhagen, Denmark 
 \endaddress 
\email thorup\@math.ku.dk
 \endemail
\catcode`\@=11
   \def\subjyear@{2010}%
\catcode`\@=\active
  \abstract\nofrills 
 Let $k$ be an uncountable field.  We prove that the polynomial ring
$R:=k[X_1,\dots,X_n]$ in $n\ge 2$ variables over $k$ is complete in
its adic topology.  In addition we prove that also the localization
$R_{\goth m}$ at a maximal ideal $\goth m\subset R$ is adically
complete. The first result settles an old conjecture of C.~U.~Jensen,
the second a conjecture of L. Gruson.  Our proofs are based on a
result of Gruson stating (in two variables) that $R_{\goth m}$ is
adically complete when $R=k[X_1,X_2]$ and $\goth m=(X_1,X_2)$.
 \endabstract
\endtopmatter

\head Introduction
\endhead

\noindent
\subhead \sectno.1\endsubhead
 Consider for a field $k$ and a given integer $n\ge 0$ the polynomial
ring $R:=k[X_1,\dots,X_n]$ in $n$ variables, and its field of
fractions $K:=k(X_1,\dots,X_n)$.  Set $d=0$ if $k$ is finite and
define $d$ by the cardinality equation $|k|=\aleph_d$ if $k$ is
infinite.  The following conjecture in its full generality was
formulated by L. Gruson (priv\.~com\., 2013).

 \proclaim{Conjecture} In the notation above,
 $\Ext_R^i(K,R)\ne0 \iff i=\inf\{d+1,n\}$. 
 \endproclaim
 The conjecture is trivially true for $n=0$ where $R=K=k$ and the
infimum equals $0$.  It is also true for $n=1$ (where $R$ is a PID.
and the infimum equals $1$; the Ext may be computed from the injective
resolution $0\to R\to K\to K/R\to 0$).

 In addition, the conjecture is trivially true if $i=0$, since the
infimum equals $0$ iff $n=0$. 

The conjecture has an obvious analogue obtained by replacing the
polynomial ring $R=k[X_1,\dots,X_n]$ by its localization $R_{\goth m}$
at a maximal ideal $\goth m$.  

\subhead\sectno.2\endsubhead
 In this note we consider the conjectures only for $i=1$. They were
formulated some 40 years ago, Conjecture \sectno.2b partly by Gruson
\cite{G, p.~254}, and Conjecture \sectno.2a by C. U. Jensen \cite{J,
p.~833}, inspired by the work of Gruson.

\proclaim{Conjectures} Let $R:=k[X_1,\dots,X_n]$ be the
polynomial ring, and $\goth m\subset R$ a maximal ideal.
Then the following bi-implications hold:
\roster
 \item"\bf 2a. "  $\Ext^1_R(K,R)\ne 0\iff n=1\text{ or } |k|\le
\aleph_0$. 
 \item"\bf 2b. " $\Ext^1_{R_{\goth m}}(K,{R_{\goth m}})\ne 0\iff 
 n=1\text{ or }  |k|\le \aleph_0$.
\endroster
\endproclaim

\subhead \sectno.3
\endsubhead
 The $\Ext$'s in the conjectures make sense for a wider class of rings,
and we fix for the rest of this paper an integral domain $R$ with
field of fractions $K$. We assume throughout that $R$ is noetherian,
and not a field; in particular, $\bigcap_{s\ne 0} sR=(0)$ and
$\Hom_R(K,R)=0$.  Let $S:=R\setminus\{0\}$ be the set of non-zero
elements of $R$, pre-ordered by divisibility: $s'\,|\,s$ iff
$sR\subseteq s'R$.  We denote by $\varprojlim^{(i)}_S$ the $i$'the
derived functor of the limit functor $\varprojlim_S$ on the category
of inverse $S$-systems of $R$-modules. 

The modules $\Ext_R^i(K,R)$ of the conjectures are related to the
$\varprojlim^{(i)}$ by well-known results, see \cite{G, p.~251--52}:
For $i\ge 2$ there are natural isomorphisms 
 $\Ext^{i}_R(K,R)\simeq \varprojlim^{(i-1)}_S R/sR$, and for
$i=1$ there is an exact sequence,
 $$0\to R@>c(R)>> \varprojlim_{s\in S}R/sR @>>> \Ext^1_R(K,R)\to 0.
\tag\sectno.3.1
 $$
 The set of principal ideals $sR$ for $s\in S$ is cofinal in the set
of all non-zero ideals of $R$.  Hence the topology defined by the
ideals $sR$ for $s\in S$ is the {\it adic topology\/} on $R$, and the
limit in \pref{\sectno.3.1} is the {\it adic completion\/} of $R$; we
denote it by $\widehat R$, and we will simply call $R$ {\it complete\/} if
the {\it canonical injection} $c(R)$ in \pref{\sectno.3.1} is an
isomorphism.  As it follows from the exact sequence
\pref{\sectno.3.1}, $R$ is complete iff $\Ext^1_R(K,R)=0$. 

Since $R$ is not a field it follows easily that the completion $\widehat
R$ is uncountable.  If the field $k$ is finite or countable (and $n\ge
1$) then the polynomial ring $R=k[X_1,\dots,X_n]$ and its localization
$R_{\goth m}$ are countable, and hence they are not complete.  In
other words, the assertions of Conjectures \sectno.2a and \sectno.2b
hold if $|k|\le \aleph_0$. As noted above, they also hold when $n\le
1$.  So the remaining cases of the conjectures are the following.

 \proclaim{Conjectures} Let $R:=k[X_1,\dots,X_n]$ be the
polynomial ring where $|k|\ge \aleph_1$ and $n\ge 2$. Then:
\roster
 \item"\bf\sectno.3a." {\rm(C. U. Jensen \cite{J, p.~833})} $R$ is
complete.
 \item"\bf\sectno.3b." The localization $R_{\goth m}$ of $R$ at any
maximal ideal $\goth m\subset R$ is complete.
 \endroster
 \endproclaim

 The main result of this paper is the verification of the two
conjectures. In fact, both conjectures are implied by a single result.

 \proclaim{Theorem \sectno.4} 
 Assume that $|k|\ge \aleph_1$, that $n\ge 2$, and that $R=U^{-1}R_0$
is a localization of $R_0=k[X_1,\dots,X_n]$ with a multiplicative
subset $U\subset R_0$. In addition, assume that every maximal ideal of
$R$ contracts to a maximal ideal of $R_0$. Then $R$ is complete.
 \endproclaim

The key ingredient in our proof is the following local result in two variables. 

 \proclaim{Proposition \sectno.5}
 {\rm(L. Gruson \cite{G, Proposition 3.2, p.~252})} Conjecture\/ {\rm
\sectno.3b} holds for $n=2$ and $\goth m=(X_1,X_2)$. 
 \endproclaim

\head Lemmas
\endhead

\subhead \sectno.6\endsubhead 
 Our argument is based on a series of lemmas, some of which are valid
in a more general context, and we keep the setup of Section \sectno.3.
First we compare for a multiplicative subset $T\subset R$ the
completions of $R$ and $T^{-1}R$. The ideals of $T^{-1}R$ generated by
the elements of $S$ form a cofinal subset of non-zero ideals.  Hence
the inclusion $R\hookrightarrow T^{-1}R$ is continuous, and 
there is an induced $R$-linear map of completions,
 $$\widehat R=\varprojlim_{s\in S} R/sR\to \varprojlim_{s\in S}
T^{-1}R/sT^{-1}R=\widehat {T^{-1}R}.\tag\sectno.6.1
 $$
  For $s\in S$ let $\goth a_s\supseteq sR$ denote the ideal of $R$
such that $\goth a_s/sR$ is the kernel of the map $R/sR\to
T^{-1}R/sT^{-1}R$.  Then the kernel of the map \pref{\sectno.6.1}
is the limit $L:=\varprojlim_{s\in S} \goth a_s/sR$.
  Clearly, for $a\in R$ we have
$a\in \goth a_s$ iff there exists an element $t\in T$ 
such that $ta\in sR$. 

 \proclaim{Lemma \sectno.7} 
 Assume that $R$ is a UFD, and let $T\subseteq S$ be a multiplicative
saturated subset. Consider the localization $R\subseteq T^{-1}R$
and the induced map of completions $\widehat R\to \widehat {T^{-1}R}$.
 Then the induced map is injective iff for every prime element
$t\in T$ there exists a prime element $p\notin T$ such that the ideal
$(t,p)R$ is proper: $(t,p)R\subset R$. 
 \endproclaim
 \demo{Proof} Recall that saturation means that any divisor of an
element of $T$ belongs to $T$ or, equivalently since $R$ is a UFD, $T$
is the submonoid of $S$ generated by a subset of prime elements.  Let
$P$ be the monoid generated by the prime elements outside $T$.
Moreover, let $T_0$ be the submonoid of $T$ consisting of elements
$t\in T$ such that $(t,p)R=R$ for all $p\in P$.  For $t\in T$ we write
$t_0$ for the largest divisor in $T_0$ of $t$, determined by a
factorization $t=t_0t'$ where $t_0\in T_0$ and $t'$ has all prime
divisors outside $T_0$.  

In this notation the Lemma asserts for the kernel $L$ of the induced
map that $L=0$ iff $T_0$ contains no prime elements. Hence the
assertion of the Lemma is a consequence of the following equation for
the kernel:
 $$L\simeq\varprojlim _{t_0\in T_0}R/t_0R.\tag\sectno.7.1
 $$
 To prove \pref{\sectno.7.1} note first that up to units in
$R^*$, the monoid $S$ is the product of $T$ and $P$, and for $s\in S$
we write $s=tp$ for the corresponding factorization into factors
$t\in T$ and $p\in P$. By unique factorization, it follows from the
description of the ideal $\goth a_s$ above that
$\goth a_s=pR$.  Consequently,
 $$\goth a_s/sR=pR/tpR\simeq R/tR.\tag\sectno.7.2
 $$
  Under the isomorphisms \pref{\sectno.7.2}, the transition map $\goth
a_s/sR\to \goth a_{s'}/s'R$ for $s'\,|\,s$ is the map $R/tR\to R/t'R$
induced by multiplication by $p/p'$. It follows from
\pref{\sectno.7.2} that
 $$L=\varprojlim_{s\in S}\goth a_s/sR
\simeq \varprojlim_{t\in T}\varprojlim_{p\in
P} R/tR.\tag\sectno.7.3
 $$
 Fix $t\in T$ and consider the inner limit in \pref{\sectno.7.3}.
 We claim that
 $$\varprojlim_{p\in P}R/tR=R/t_0R.\tag\sectno.7.4
 $$
 The transition maps for the limit in \pref{\sectno.7.4} are
multiplications by elements $p\in P$ on the $R$-module $R/tR$.  By
unique factorization, the multiplications are injective.  Therefore,
the limit is the intersection of the images of the multiplications. 

Clearly, if $t\in T_0$ then the multiplications are bijective; hence
the intersection is equal to $R/tR$, and \pref{\sectno.7.4} holds
since $t=t_0$.  Assume next that $t$ is a prime element outside $T_0$.
Then multiplication by some $p\in P$ has an image contained in a
proper ideal of $R/tR$. Hence the intersection of the images is
contained in the intersection of the powers of a proper ideal of
$R/tR$. Since $R/tR$ is an integral domain, the intersection equals
$0$ by Krull's Intersection Theorem, and hence \pref{\sectno.7.4}
holds since $t_0=1$.

In general, we factorize $t=t_0t'$ where $t'$ has all prime divisors
outside $T_0$, and use the exact sequence $0\to R/t'R\to R/tR\to
R/t_0R\to 0$.  From the previous considerations it follows first that
the intersection of the images on $R/t'R$ is equal to $0$, and next
that the intersection of the images on $R/tR$ maps isomorphically onto
$R/t_0R$.  Hence \pref{\sectno.7.4} holds in general.

Clearly \pref{\sectno.7.4} and \pref{\sectno.7.3} imply
\pref{\sectno.7.1}.
 \enddemo

 \proclaim{Lemma \sectno.8} 
 Assume for every maximal ideal $\goth m$ of $R$ that the induced map
$\widehat R\to \widehat {R_{\goth m}}$ is injective and that $R_{\goth m}$
is complete. Then $R$ is complete. 
 \endproclaim
\demo{Proof}
 By the second assumption, $R_{\goth m}=\widehat {R_{\goth m}}$.
Hence, by the first assumption, $\widehat R$ embeds into $\bigcap
R_{\goth m}=R$. Thus $R=\widehat R$. 
\enddemo

 \proclaim{Lemma \sectno.9} Assume that $R\subseteq R'$ is a subring
of an integral domain $R'$ such that $R'$ is integral over $R$ and
free as an $R$-module.  Assume that $R'$ is complete.  Then $R$ is
complete.
 \endproclaim
 \demo{Proof}
 Every non-zero ideal of $R'$ contracts to a non-zero ideal of $R$
since $R'$ is integral over $R$. In other words, the inclusion $R\to
R'$ is continuous. Hence there is an induced map of completions
$\widehat {R}\to \widehat{R^{\prime\,}}$, and an induced $R'$-linear map
$R'\otimes_R\widehat R\to \widehat{R^{\prime\,}}$.

 We have to prove that the canonical injection $c=c(R)\:R\to
\widehat R$ is an isomorphism. Since $R'$ is free over $R$ it suffices
to prove that the map $R'\otimes_{R}c\:R'\to R'\otimes_{R}
\widehat{R}$ is an isomorphism. Clearly, the canonical injection
$c(R')\:R'\to \widehat{R'} $ factors:
 $$\CD R' @=  R'\\
@VR'\otimes_{R}cVV @VVc(R')V\\
R'\otimes_{R} \widehat{R}@>>> \widehat{R'}\rlap.
\endCD
 $$
 The bottom map is the canonical map $R'\otimes_{R}\varprojlim V_s\to
\varprojlim (R'\otimes_{R}V_s)$ defined for any inverse $S$-system of
$R$-modules $(V_s)$.  It is injective, since $R'$ is free over $R$. The right
vertical map is an isomorphism by assumption. Therefore, $R'\otimes
_{R}c$ is an isomorphism. 
 \enddemo

\proclaim{Lemma \sectno.10} 
 Assume that $R$ is a localization of
$R_0=k[X_1,\dots,X_n]$ such that every maximal ideal
of $R$ contracts to a maximal ideal of $R_0$. Let $\goth p\subset R$
be a prime ideal of height at least $2$. Then the induced map of
completions $\widehat R\to \widehat{R_{\goth p}}$ is injective.
 \endproclaim
 \demo{Proof} Indeed, as is well-known, the localization $R$ is a UFD:
its prime elements are, up to units in $R^*$, those irreducible
polynomials in $R_0$ that are non-units of $R$.  To apply Lemma
\cref{\sectno.7}, let $t$ be a prime element in $R\setminus \goth p$.
We have to prove that there exists a prime element in $\goth p$ such
that the ideal
$(t,p)R$ is proper. Take any maximal ideal $\goth m\subset R$ with
$t\in \goth m$. Apply the following Sublemma to the contractions
$\goth m_0=R_0\cap \goth m$ and $\goth p_0=R_0\cap \goth p$. It
follows that there exists an irreducible polynomial $p$ in 
$\goth m_0\cap \goth p_0$. Then $p$ is a prime element in $\goth p$, 
and $(t,p)R$ is a proper ideal, since $(t,p)R\subseteq \goth m$. 
 \enddemo

 \proclaim{Sublemma} Let $R=k[X_1,\dots,X_n]$ be the
polynomial ring, let $\goth p\subset R$ be a prime ideal of height
$h\ge 2$, and let $\goth m\subset R$ be a maximal ideal.  Then the
intersection $\goth p\cap \goth m$ contains a prime ideal $\goth q$ of
height $h-1$. In particular, $\goth m\cap\goth p$ contains an
irreducible polynomial. 
 \endproclaim
 \demo{Proof} The assertion is trivial if $\goth p\subseteq \goth m$
so we may assume that $\goth p\not\subseteq \goth m$.  Assume first
that $k$ is algebraically closed.  Then $\goth p$ is the ideal of an
irreducible variety $V$, and $\goth m$ is the ideal of a point $q$. By
assumption $q\notin V$. Hence the linear join of $q$ and $V$ (the cone
with base $V$ and vertex $q$) is an irreducible subvariety $W$ of
dimension equal to $\dim V+1$.  Therefore, the ideal $\goth q$ of $W$
is a prime ideal with the required properties.

The general case is reduced to the previous case as follows: Consider
the embedding $R\hookrightarrow \bar R$ where $\bar R$ is the
polynomial ring over the algebraic closure of $k$.  The embedding is
integral, and $R$ is a UFD.  Hence, by the usual dimension theory for
polynomial rings, $\goth p$ and $\goth m$ are contractions of prime
ideals $\overline{\goth p}$ and $\overline{\goth m}$ of $\bar R$; if
$\overline{\goth q}\subseteq \overline{\goth p}\cap \overline{\goth
m}$ is a prime of height $h-1$, then the contraction $\goth q:=R\cap
\overline{\goth q}$ has the required property.
 \enddemo

\remark{\bf Note \sectno.11} (1)
 The proof of Lemma \cref{\sectno.10} is particularly simple in the
special case: $R=k[X_1,X_2]$, $k$ is algebraically closed, and
$\goth p=(X_1,X_2)$. Indeed, for an irreducible polynomial $t$ outside $\goth p$
take a zero $\alpha=(\alpha _1,\alpha _2)$ of $t$ and take $p:=\alpha
_2X_1-\alpha _1X_2$. Then $t$ and $p$ belong to the maximal ideal
$\goth m_\alpha=(X_1-\alpha_1,X_2-\alpha_2)$, and $p$ is irreducible
since $\alpha\ne (0,0)$. The special case is sufficient
for a proof of Conjecture \cref{3a} alone, see Note
\cref{\sectno.14(2)}.

(2) It is also worthwhile to note that the conclusion in Lemma \sectno.10
is wrong for prime ideals $\goth p$ of height $1$: For
$R=k[X_1,\dots,X_n]$ and a prime ideal $\goth p$ of height $1$, say
$\goth p=pR$, the induced map $\widehat R\to \widehat{R_{\goth p}}$ is
not injective. Indeed, the polynomial $p+1$ is not a constant, and
hence for any irreducible divisor $t$ in $1+p$ we have $t\notin\goth
p$ and $(t,p)R=R$. Hence, by Lemma \cref{\sectno.7}, the map is not
injective. 
 \endremark

\head Proofs of the main results
\endhead

 \proclaim{Lemma \sectno.12} 
 Let $R:=k[X_1,X_2]$ be the polynomial ring in two variables where
$|k|\ge \aleph_1$, and let $\goth m$ be any maximal ideal of $R$.
Then $R$ and $R_{\goth m}$ are complete.
 \endproclaim
 \demo{Proof} The second assertion is a generalization of Gruson's
local result.  First, if $k$ is algebraically closed, then $R_{\goth
m}$ is complete. Indeed, then $\goth m=(X_1-\alpha _1,X_2-\alpha _2)$
with $\alpha _1,\alpha _2\in k$, and the completeness of $R_{\goth m}$
follows from the local result (Proposition \cref{\sectno.5}) by a
change of coordinates.  

To prove the results in general, embed $k$ in the algebraic closure
$\bar k$. Let $\bar R:=\bar k[X_1,X_2]$, let $R':=R_{\goth m}$ and
$\bar R':=\bar R_{\goth m}$.  With $U:=R\setminus \goth m$ we have
$R'=U^{-1}R$ and $\bar R'=U^{-1}\bar R$. The maximal ideals of $\bar
R'$ are the ideals generated by maximal ideals $\bar {\goth m}\subset
\bar R$ lying over $\goth m$. Moreover, the localization of $\bar R'$
at the maximal ideal $\bar{\goth m}\bar R'$ is equal to $\bar
R_{\bar{\goth m}}$, and hence complete by the first case.  In
addition, the map of completions induced by $\bar R'\hookrightarrow
\bar R_{\bar{\goth m}}$ is injective by Lemma \cref{\sectno.10}.
 Therefore, by Lemma \cref{\sectno.8}, the ring $\bar R'$ is complete.
Finally, $\bar R'=\bar R\otimes_RR'$ is integral and free over $R'$.
Hence, by Lemma \cref{\sectno.9}, $R'$ is complete.  Similarly, since
$\bar R_{\bar{\goth m}}$ is complete for all maximal ideals
$\bar{\goth m}$ of $\bar R$, it follows first that $\bar R$ is
complete, and next the $R$ is complete.

 \enddemo

 \proclaim{Theorem \sectno.13}
 Assume that $|k|\ge \aleph_1$, that $n\ge 2$, and that $R=U^{-1}R_0$
is a localization of $R_0=k[X_1,\dots,X_n]$ with a multiplicative
subset $U\subset R_0$. In addition, assume that every maximal ideal of $R$
contracts to a maximal ideal of $R_0$. Then $R$ is complete.
 \endproclaim \demo{Proof} Clearly $R$ is a UFD, and hence equal to
the intersection of the localizations $R_{\goth q}$ over all prime
ideals $\goth q$ of height $1$. Moreover, every height $1$ prime ideal
is contained in a height $2$ prime ideal, since $R$ is catenary and
all maximal ideals have height $n\ge 2$. Therefore, $R$ is the
intersection over all prime ideals $\goth p$ of height $2$:
 $$R=\bigcap_{\operatorname{ht}\goth p=2}R_{\goth p}.\tag\sectno.13.1
 $$
 For every prime ideal $\goth p$ of height $2$ it follows from Lemma
\cref{\sectno.10} that the induced map of completions $\widehat R\to
\widehat{R_{\goth p}}$ is injective.  Therefore, by Equation
\pref{\sectno.13.1}, to prove that $R=\widehat R$, it suffices to
prove for every height $2$ prime ideal $\goth p$ of $R$ that $R_{\goth
p}$ is complete. Clearly, the latter completeness follows from Lemma
\cref{\sectno.12} using the following standard observation on
localizations for $h=2$: {\sl If $R$ is a localization of
$R_0=k[X_1,\dots,X_n]$ then any localization $R_{\goth p}$ at a prime
ideal $\goth p\subset R$ of height $h\ge 1$ may be obtained, after a
renumbering of the variables, as the localization at a maximal ideal
of the polynomial ring,}
 $$k(X_{h+1},\dots,X_n)[X_1,\dots,X_h].\tag\sectno.13.2
 $$
 To justify the observation, note first that the prime ideal $\goth
p\subset R$ is generated by a prime ideal $\goth p_0\subset R_0$, and
$R_{\goth p}=(R_0)_{\goth p_0}$.
 Hence we may assume that $R=k[X_1,\dots,X_n]$. The quotient $R/\goth
p$ has transcendence degree $n-h$ over $k$ since $\goth p$ has height
$h$.  Consequently there are $n-h$ among the variables, say
$X_{h+1},\dots,X_{n}$, whose classes modulo $\goth p$ are
algebraically independent, or equivalently, such that
$k[X_{h+1},\dots,X_n]\cap \goth p=(0)$. Localization of $R$ with the
monoid of non-zero polynomials in $X_{h+1},\dots,X_n$ yields the ring
\pref{\sectno.13.2}, and so $R_{\goth p}$ may be obtained by
localization of \pref{\sectno.13.2} at the ideal generated by $\goth
p$.  The latter ideal is a prime ideal of height $h$, and hence a
maximal ideal.  Thus the observation has been justified.
 \enddemo

\remark{\bf Note \sectno.14} 
 (1) Clearly Theorem \cref{\sectno.13} implies the two conjectures 3a and
3b.  In addition, it follows from the observation at the end of the
previous proof that Conjecture 3b implies, when $|k|\ge \aleph_1$,
that the localization $R_{\goth p}$ of $R:=k[X_1,\dots,X_n]$ at any
prime ideal $\goth p$ of height $h\ge 2$ is complete. In particular,
the rings in Theorem \cref{\sectno.13} do not exhaust the list of
complete subrings of $k[X_1,\dots,X_n]$.

 (2) For a proof of Conjecture \cref{\sectno.3a} alone, the arguments
can be simplified.  First, the proof of Lemma \cref{\sectno.12} for
$R=k[X_1,X_2]$ uses only the special case of Lemma \cref{\sectno.10}
mentioned in Note \cref{\sectno.11}(1). Next, for $R=k[X_1,\dots,X_n]$
with $|k|\ge \aleph_1$ and $n\ge 3$ a direct proof of completeness is
the following:

Denote by $T_{12}\subset S$ the multiplicative subset of polynomials containing
neither $X_1$ nor $X_2$, that is, $T_{12}$ is the set of non-zero
polynomials in $k[X_3,\dots,X_n]$. Then the localization,
 $$T_{12}^{-1}R=k(X_3,\dots,X_n)[X_1,X_2],
 $$
 is complete by Lemma \cref{\sectno.12}. Moreover, it follows
immediately from Lemma \cref{\sectno.7} that the inclusion
$R\hookrightarrow T_{12}^{-1}R$ induces an injection on the
completions; indeed, for any irreducible $t\in T_{12}$ take $p=X_1$.
Similarly, with an obvious notation we obtain for any $i=3,\dots,n$ an
inclusion $\widehat R\hookrightarrow T_{1i}^{-1}R$, and hence an
inclusion,
 $$\widehat R\hookrightarrow \bigcap_{i=2}^n T_{1i}^{-1}R.
 $$
 Obviously, the intersection on the right side equals $R$. Thus
$\widehat R=R$. 
 \endremark

\widestnumber\key M

\enddocument